\documentclass{amsart}
\title[On metric spaces with the properties of de Groot and Nagata]{On metric spaces with the properties of de Groot and Nagata in dimension one}
\author[T.Banakh, D.Repov\v s and I.Zarichnyi]{Taras Banakh, Du\v san Repov\v s and Ihor Zarichnyi}
\address{T.Banakh: Department of Mathematics, Ivan Franko National University of Lviv, Ukraine; and Instytut Matematyki, Uniwersytet Humanistyczno Przyrodniczy im. Jana Kochanowskiego w Kielcach, Poland}
\email{tbanakh@yahoo.com}
\address{D.Repov\v s: Faculty of Mathematics and Physics,
and Faculty of Education, University of Ljubljana, P.O.B.2964, Ljubljana, Slovenia}
\email{dusan.repovs@guest.arnes.si}
\address{I.Zarichnyi: Department of Mathematics, Ivan Franko National University of Lviv, Ukraine}
\email{ihor.zarichnyj@gmail.com}
\thanks{This research was supported by Slovenian Research Agency grant P1-0292-0101, J1-9643-0101 and BI-UA/07-08-001.}

\date{\today}

\newcommand{\e}{\varepsilon}
\newcommand{\GP[1]}{\mathsf{GP}_{#1}}
\newcommand{\NP[1]}{\mathsf{NP}_{#1}}
\newcommand{\IR}{\mathbb R}
\newcommand{\IQ}{\mathbb Q}
\newcommand{\II}{\mathbb I}
\newcommand{\IC}{\mathbb C}
\newcommand{\w}{\omega}
\newcommand{\dist}{d}
\newcommand{\Lip}{\mathrm{Lip}}
\newcommand{\dens}{\mathrm{dens}}
\newcommand{\F}{\mathcal F}

\newcommand{\diam}{\mathrm{diam}}
\newcommand{\Dist}{\mathrm{Dist}}

\newcommand{\sign}{\mathrm{sign}}

\newtheorem{theorem}{Theorem}[section]
\newtheorem{exth}[theorem]{Extension Theorem}
\newtheorem{univth}[theorem]{Universality Theorem}
\newtheorem{corollary}[theorem]{Corollary}
\newtheorem{lemma}[theorem]{Lemma}
\newtheorem{proposition}[theorem]{Proposition}
\theoremstyle{definition}
\newtheorem{problem}[theorem]{Problem}
\newtheorem{example}[theorem]{Example}

\begin{document}
\begin{abstract} A metric space $(X,d)$ has the de Groot property $\GP[n]$ if for any points $x_0,x_1,\dots,x_{n+2}\in X$ there are positive indices $i,j,k\le n+2$ such that $i\ne j$ and $d(x_i,x_j)\le d(x_0,x_k)$. If, in addition, $k\in\{i,j\}$ then $X$ is said to have the Nagata property $\NP[n]$. It is known that a compact metrizable space $X$ has dimension $\dim(X)\le n$ iff $X$ has an admissible $\GP[n]$-metric iff $X$ has an admissible $\NP[n]$-metric.

We prove that an embedding $f:(0,1)\to X$ of the interval $(0,1)\subset \IR$ into a locally connected metric space $X$ with property $\GP[1]$ (resp. $\NP[1]$) is open,  provided $f$ is an isometric embedding (resp. $f$ has distortion $\Dist(f)=\|f\|_\Lip\cdot\|f^{-1}\|_\Lip<2$).  This implies that the Euclidean metric cannot be extended from the interval $[-1,1]$ to an admissible $\GP[1]$-metric on the triode $T=[-1,1]\cup[0,i]$. Another corollary says that a topologically homogeneous $\GP[1]$-space cannot contain an isometric copy of the interval $(0,1)$ and a topological copy of the triode $T$ simultaneously. Also we prove that a $\GP[1]$-metric space $X$ containing an isometric copy of each compact $\NP[1]$-metric space has density $\ge \mathfrak c$.
\end{abstract}

\maketitle

\section{Introduction}
In this paper we shall be interested in structural properties of metric spaces possessing the properties introduced by J. de Groot \cite{dG} and J.~Nagata \cite{Nag}. 

Let $n$ be a non-negative integer. A metric $d$ on $X$ is said to have the {\em de Groot property} $\GP[n]$ if for any $n+3$ points $x_0,x_1,\dots,x_{n+2}\in X$ there is a triplet of indices $i,j,k\in\{1,\dots,n+2\}$ such that
$$d(x_i,x_j)\le d(x_0,x_k)\mbox{ \ and \ $i\ne j$}.
$$
If, in addition, $k\in\{i,j\}$, then we say that the metric $d$ has the {\em Nagata property} $\NP[n]$ or that $d$ is an $\NP[n]$-metric. It is clear that each $\NP[n]$-metric is also a $\GP[n]$-metric. In the Engelking's monograph \cite{Eng} the properties of Nagata and de Groot are denoted by $(\mu_4)$ and $(\mu_5')$, respectively.
Those properties also are discussed in the Nagata's book \cite[V.3]{NagB}.

According to \cite{dG} and \cite{Nag}, for a separable metrizable space $X$ the following conditions are equivalent:
\begin{itemize}
\item $X$ has the covering dimension $\dim(X)\le n$;
\item the topology of $X$ is generated by an $\NP[n]$-metric on $X$;
\item the topology of $X$ is generated by a totally bounded $\GP[n]$-metric on $X$.
\end{itemize}
 In fact, the equivalence of the first two conditions hold for any metrizable space $X$. On the other hand, it is an open problem due to de Groot \cite{dG} if the existence of an admissible $\GP[n]$-metric on a (separable) space $X$ implies $\dim(X)\le n$, see \cite[p.231]{Eng}. We recall that a metric $d$ on a topological space $X$ is said to be {\em admissible} if it generates the topology of $X$.

By \cite[4.2.D]{Eng}, a metric $d$ has the $\GP[0]$-property if and only if it has the $\NP[0]$-property if and only if the metric $d$ satisfies the strong triangle inequality $$d(x_1,x_2)\le\max\{d(x_0,x_1),d(x_0,x_2)\}$$ for all points $x_0,x_1,x_2\in X$. The latter means that $d$ is an {\em ultrametric}. Thus both $\NP[n]$-metric and $\GP[n]$-metric are higher dimensional analogs of ultrametric.

Due to efforts of many mathematicians the structure of ultrametric spaces is quite well understood. We shall recall  two results: an Extension Theorem and a Universality Theorem.

\begin{exth}\label{exth} Each admissible ultrametric defined on a closed subspace $A$ of an zero-dimensional compact metrizable space $X$ extends to an admissible ultrametric on $X$.
\end{exth}

This theorem follows from its uniform version proved by Ellis in \cite{Ellis} or its ``simultaneous'' version proved by Tymchatyn and Zarichnyi \cite{TZ}. The other theorem is due to  A.Liman and V.Liman \cite{LL} and concerns universal ultrametric spaces. We define a (topological) metric space $X$ to be ({\em topologically}) {\em homogeneous} if for any two points $x,y\in X$ there is an isometry (a  homeomorphism) $h:X\to X$ such that $h(x)=y$.

\begin{univth}\label{univtheorem} For each cardinal $\kappa$ there is a (homogeneous) ultrametric space $LM_\kappa$ of weight $\kappa^\w$ containing an isometric copy of each ultrametric space of weight $\le\kappa$.
\end{univth}

The universal space $LM_\kappa$ in Theorem~\ref{univtheorem} can be constructed as follows: take any Abelian group $G$ of size $|G|=\kappa$, let $\IQ_+$ be the set of all positive rational numbers, and let $LM_\kappa$ be the space of all maps $f:\IQ_+\to G$ which are eventually zero, in the sense that $f(x)$ is zero for all sufficiently large rational numbers $x\in\IQ_+$. The space $LM_\kappa$ endowed with the ultrametric $d(f,g)=\sup\{x\in\IQ_+:f(x)\ne g(x)\}$ (where $\sup\emptyset=0$) has the structure of an Abelian group and therefore is metrically homogeneous.
\smallskip

It is natural to ask if these two theorems have analogues for $\GP[n]$ or $\NP[n]$-metrics. As we shall see later, the answer is negative already for $n=1$.
To construct a suitable counterexample we shall first study the structure of $\GP[1]$-spaces $X$ in a neighborhood of an isometrically embedded interval $(0,1)\subset X$.
\smallskip

\begin{theorem}\label{c1} If a $\GP[1]$-metric space $X$ is locally connected, then each subset $I\subset X$, isometric to an interval $(a,b)\subset\IR$, is open in $X$.
\end{theorem}

This theorem will be proved in Section~\ref{s2}. Now we discuss some of its corollaries.

By the {\em triode\/} we understand the subspace $$T=[-1,1]\cup[0,i]$$
 of the complex plane $\IC$. By Nagata's Theorem \cite{Nag}, the triode $T$ carries an admissible $\NP[1]$-metric. Nonetheless, such a metric cannot restrict to the Euclidean metric on the interval $[-1,1]\subset T$ because the interval $(-1,1)$ is not open in the triode. Thus we obtain:

\begin{corollary}\label{c2} The Euclidean metric on the interval $[-1,1]$ has the Nagata property $\NP[1]$ but cannot be extended to an admissible $\GP[1]$-metric on the triode $T$.
\end{corollary}

Therefore, Extension Theorem~\ref{exth} cannot be generalized to metric spaces with the property $\NP[n]$ or $\GP[n]$ for $n\ge 1$.
Next, we show that the same concerns Universality Theorem~\ref{univtheorem}: its homogeneous version cannot be generalized to higher dimensions.

\begin{corollary}\label{c3} If a $\GP[1]$-metric space $X$ contains both an isometric copy of the interval $[0,1]$ and a topological copy of the triode $T$, then $X$ is not topologically homogeneous.
\end{corollary}

\begin{proof} Let $[0,1]\subset X$ be an isometric copy of the interval $[0,1]$.  Assuming that  $X$ is topologically homogeneous and $X$ contains a topological copy of the triode $T$, we can find a topological embedding $f:T\to X$ such that $f(0)=\frac12\in[0,1]\subset X$. Since the triode does not embed into the interval $[0,1]$, the point $1/2$ is not an interior point of the interval $(0,1)$ in the locally connected subspace $Y=[0,1]\cup f(T)$ of the $\GP[1]$-space $X$.  This contradicts Theorem~\ref{c1}.
\end{proof}

In spite of the negative result in Corollary~\ref{c3}, we do not know the answer to the following

\begin{problem}\label{pr1.6} Is it true that for each infinite cardinal $\kappa$ there is a $\GP[1]$-metric space $U$ of weight $\kappa^\w$ that contains an isometric copy of each $\NP[1]$-metric space $X$ of weight $\le\kappa$?
\end{problem}

The weight $\kappa^\w$ in Problem~\ref{pr1.6} cannot be replaced by $\kappa$ because of the following theorem that will be proved in Section~\ref{s3}.

\begin{theorem}\label{dens} If a $\GP[1]$-metric space $X$ contains an isometric copy of each compact $\NP[1]$-metric space, then $X$  has density $\dens(X)\ge\mathfrak c$.
\end{theorem}

Now let us return to Theorem~\ref{c1}. It implies that no non-open arc $I$ in a locally connected $\GP[1]$-metric space $(X,d)$
is  isometric to an interval $(a,b)\subset \IR$.  We can ask
how much the metric $d$ restricted to $I$ differs from the Euclidean
metric on $I$. We can measure this distance using the notion of
the distortion.

By the {\em distortion} of an injective map $f:X\to Y$ between metric
spaces $(X,d_X)$ and $(Y,d_Y)$ we understand the
(finite or infinite) number
$$\Dist(f)=\|f\|_\Lip\cdot\|f^{-1}\|_\Lip$$where
$$\|f\|_{\Lip}=\sup_{x\ne x'}\frac{d_Y(f(x),f(x')}{d_X(x,x')}$$
is the Lipschitz constant of $f$ (if $|X|\le 1$, then $\|f\|_\Lip$
is not defined, so we put $\Dist(f)=1$).
The notion of distortion is widely used in studying the
embeddability problems of metric spaces, see \cite{naor},
\cite{mat1}, \cite{mat2}, \cite{Men}.

It can be shown that an embedding $f:X\to Y$ of a  metric space $X$ into a metric space $Y$ has distortion $\Dist(f)=1$ if and only if $f$ is a {\em similarity}, which means that $d_Y(f(x),f(x'))=\|f\|_\Lip\cdot d_X(x,x')$ for all $x,x'\in X$.

In terms on the distortion, Theorem~\ref{c1} can be written as follows.

\begin{corollary}\label{t5} Let $X$ be a locally connected metric space with property $\GP[1]$. Each embedding $f:(0,1)\to X$ with distortion $\Dist(f)=1$ is open.
\end{corollary}

\begin{proof} Let $f:(0,1)\to X$ be an embedding with distortion $\Dist(f)=1$. Let $C=\|f\|_\Lip$ and $$g:(0,C)\to(0,1),\;\;g:t\mapsto t/C,$$ be the similarity mapping having the Lipschitz constant $\|g\|_\Lip=1/C$. It follows that the composition $f\circ g:(0,C)\to X$ has distortion $$1=\Dist(f\circ g)=\|f\circ g\|_\Lip\cdot\|(f\circ g)^{-1}\|_\Lip=1.$$Since $\|f\circ g\|_\Lip=1$, we conclude that $\|(f\circ g)^{-1}\|_\Lip=1$ and hence $f\circ g:(0,C)\to X$ is an isometric embedding. By Theorem~\ref{c1}, the image $f\circ g\big((0,C)\big)=f\big((0,1)\big)$ is open in $X$.
\end{proof}

\begin{problem} Can the equality $\Dist(f)=1$ in Corollary~\ref{t5} be replaced by the inequality $\Dist(f)<2$.
\end{problem}

This problem has an affirmative solution for metric spaces with the Nagata property $\NP[1]$. The following theorem can be easily derived from Proposition~\ref{BM} and Corollary~\ref{c8} proved at the end of the paper.

\begin{theorem}\label{t6} Let $X$ be a locally connected metric space with property $\NP[1]$. Each embedding $f:(0,1)\to X$ with distortion $\Dist(f)<2$ is open.
\end{theorem}

The inequality $\Dist(f)<2$ in this theorem is best possible because of the following 
simple example.

\begin{example} On the triode $T=[-1,1]\cup[0,i]$ consider the $\NP[1]$-metric $$\rho(z,z')=\begin{cases}
|z-z'|&\mbox{if $\sign(\Re(z))=\sign(\Re(z'))$,}\\
\max\{|\Re(z)|,|\Re(z')|,\Im(z),\Im(z')\}&\mbox{otherwise}.
\end{cases}
$$It is easy to check that the identity embedding $f:[-1,1]\to (T,\rho)$ has distortion
$\Dist(f)=2$ but is not open.
\end{example}

In spite of Corollary~\ref{c2} there is a hope that the following problem (related to an approximative extension of $\NP[1]$-metrics) has an affirmative solution.

 \begin{problem} Let $A$ be a closed subspace of a 1-dimensional space $X$. Is it true that for any a admissible $\NP[1]$-metric $d_A$ on $A$ there is an admissible $\NP[1]$-metric $d_X$ on $X$ such that the identity embedding $f:(A,d_A)\to (X,d_X)$ has distortion $\Dist(f)\le 2$?
\end{problem}

\section{Isometric arcs in $\GP[1]$-metric spaces}\label{s2}

In this section we shall prove Theorem~\ref{c1}. A map $f:X\to Y$ between metric spaces is called {\em non-expanding} if its Lipschitz constant $\|f\|_{\Lip}\le 1$.
For a point $x$ of a metric space $(X,d)$ and a subset $A\subset X$ we put $d(x,A)=\inf_{a\in A}d(x,a)$.

\begin{lemma}\label{slice} Let $(X,d)$ be a $\GP[1]$-metric space containing an isometric copy of the closed interval $[0,1]$ and let \ $V=\{x\in X:\dist(x,[0,1])<\frac13\dist(x,\{0,1\})\}.$
\begin{enumerate}
\item There is a non-expanding retraction $r:V\to (0,1)$ such that $$d(x,t)=\max\{|t-r(x)|,\dist(x,[0,1])\}\mbox{ for any }x\in V,\; t\in(0,1).$$
\item For any points $x,y\in V$ with $\dist(x,[0,1])\ne\dist(y,[0,1])$ we get $$d(x,y)\ge \max\{\dist(x,[0,1]),\dist(y,[0,1])\}.$$
\end{enumerate}
\end{lemma}

\begin{proof} 1. Given any $x\in V$, let $D=\dist(x,[0,1])$ and  consider the compact subset $D(x)=\{t\in[0,1]:d(x,t)=D\}$.  We claim that $D(x)$ is a closed subinterval of $(0,1)$ of length $2D$. Let $a=\min D(x)$ and $b=\max D(x)$.

The triangle inequality implies that $d(a,b)\le d(a,x)+d(x,b)\le 2D$. It follows from $D<\frac13\dist(x,\{0,1\})$ that $d(0,a)\ge d(0,x)-d(x,a)>3D-D>D$ and similarly, $d(b,1)>D$. Let us show that $[a,a+D]\subset D(x)$. Assuming the converse, we could find a point $x_1\in (a,a+D]\setminus D(x)$.
Then for the points $$x_0=a,\;\;x_1,\;\; x_2=x, \mbox{ and }x_3=a-D$$ we would get $$
\begin{gathered}
d(x_1,x_2)>D,\; d(x_1,x_3)=D+(x_1-a)>D,\;d(x_2,x_3)>D \mbox{ and }\\
d(x_0,x_3)=d(a,a-D)=D,\; d(x_0,x_2)=d(a,x)=D,\; d(x_0,x_1)=d(a,x_1)\le D,
\end{gathered}$$
which contradicts the $\GP[1]$-property of the metric $d$.

Thus $[a,a+D]\subset D(x)$. By analogy we can prove that $[b-D,b]\subset D(x)$. Combined with $b-a\le 2D$, this implies that $[a,b]=[a,a+D]\cup[b-D,b]=D(x)$. Assuming that $b-a<2D$, we could take
$x_0$ be the midpoint of the interval $[a,b]$ and put $x_1=x$, $x_2=x_0-D$, $x_3=x_0+D$. Then
$$
\min\{d(x_1,x_2),d(x_1,x_3),d(x_2,x_3)\}>D=\max\{d(x_0,x_1),d(x_0,x_2),d(x_0,x_3)\},
$$which contradicts the $\GP[1]$-property of the metric $d$.

Therefore, $D(x)$ is a closed interval of length $2D$. Let $r(x)$ be the midpoint of this interval. Let us show that $d(x,t)=\max\{|t-r(x)|,D\}$ for all $t\in[0,1]$. This is obvious if $t\in D(x)=[a,b]$. So assume that $t\notin D(x)$. If $t<a$, then $d(t,x)\le d(t,a)+d(a,x)\le a-t+D=r(x)-t$.
On the other hand, $b-t=d(t,b)\le d(t,x)+d(x,b)=d(t,x)+D$ implies $d(t,x)\ge b-t-D=r(x)-t$. Therefore $d(x,t)=r(x)-t=\max\{|r(x)-t|,D\}$.
The case $t>b$ can be treated by analogy.

Finally, we show that the map $r:V\to(0,1)$, $r:x\mapsto r(x)$ is a non-expanding retraction. It is clear that $r(t)=t$ for any $t\in (0,1)$. Take any two points $x,y\in V$. Without loss of generality, $r(y)\ge r(x)$. Let $D_x=\dist(x,[0,1])$ and $D_y=\dist(y,[0,1])$. For the point $t=r(x)-D_x=\min D(x)$ let us observe that
\begin{multline*}
r(y)-r(x)+D_x=r(y)-t\le \max\{|r(y)-t|,D_y\}=\\
d(t,y)\le d(t,x)+d(x,y)=D_x+d(x,y)
\end{multline*} and hence $|r(y)-r(x)|=r(y)-r(x)\le d(x,y)$.
\medskip

2. Take any two points $x,y\in V$ with $D_x=\dist(x,[0,1])\ne\dist(y,[0,1])=D_y$. We need to prove that $d(x,y)\ge\max\{D_x,D_y\}$. Without loss of generality, $D_x<D_y$. Assume conversely that $d(x,y)<\max\{D_x,D_y\}=D_y$.
Observe that $$d(r(x),0)\ge d(x,0)-D_x\ge d(y,0)-d(x,y)-D_x>d(y,0)-2D_y>3D_y-2D_y=D_y$$ and hence for any real $a$ with $\max\{D_x,d(x,y)\}<a<D_y$ the point $x_1=r(x)-a\in(0,1)$ is well-defined. By analogy we can prove that $x_2=r(x)+a\in(0,1)$ is well-defined.

So we can consider the 4 points: $x_0=x$, $x_1=r(x)-a$, $x_2=r(x)+a$, $x_3=y$, and derive a contradiction with the $\GP[1]$-property of the metric $d$ because:
$$
\begin{aligned}&\min\{d(x_1,x_2),d(x_1,x_3),d(x_2,x_3)\}\ge \min\{2a,D_y,D_y\}>\\
&> \max\{a,a,d(x,y)\}\ge\max\{d(x_0,x_1),d(x_0,x_2),d(x_0,x_3)\}.
\end{aligned}$$
\end{proof}

\begin{proof}[Proof of Theorem~\ref{c1}.] Let $X$ be locally connected $\GP[1]$-metric space and $I\subset X$  a subset isometric to an open interval $(a,b)\subset\IR$.
 We need to check that each point $x_0\in I$ is an interior point of $I$ in $X$. For a sufficiently small $\e>0$ we can find an isometry $f:[0,2\e]\to I\subset X$ such that $f(\e)=x_0$. Scaling the $\GP[1]$-metric $d$ of $X$ by a suitable constant, we can assume that $\e=\frac12$. We shall identify the interval $[0,1]$ with a subinterval of $I$ and $1/2$ with the point $x_0$. Consider the neighborhood $$V=\{x\in X:\dist(x,[0,1])<\dist(x,\{0,1\})/3\}$$ of $(0,1)$ in $X$. By the local connectedness of $X$ at $x_0$,  find a connected neighborhood $C(x_0)\subset V$ of the point $x_0=1/2$. We claim that $C(x_0)\subset I$. Otherwise there would exist a point $x_1\in C(x_0)\setminus I$. Lemma~\ref{slice}(2) guarantees that the subset $$D=\{x\in C(x_0):\dist(x,[0,1])=\dist(x_1,[0,1])\}$$ is open-and-closed in $C(x_0)$, which implies that the neighborhood $C(x_0)$ is not connected and this is a contradiction.
\end{proof}

\section{Universal $\GP[1]$-spaces}\label{s3}

In this section we study universal $\GP[1]$-spaces and prove Lemma~\ref{univ} which implies Theorem~\ref{dens} announced in the Introduction.

We shall need the following (probably known)

\begin{lemma}\label{P1+0} Let $(X,d_X)$ be a $\NP[1]$-metric space and $(Y,d_Y)$ be an $\NP[0]$-metric space. Then the $\max$-metric
$$\dist\big((x,y),(x',y')\big)=\max\{d_X(x,x'),d_Y(y,y')\}$$ on the product $X\times Y$ has the Nagata property $\NP[1]$.
\end{lemma}

\begin{proof} Given any 4 points $(x_0,y_0),(x_1,y_1),(x_2,y_2),(x_3,y_3)\in X\times Y$, we need to find two distinct indices $i,j\in\{1,2,3\}$ such that
$$d\big((x_i,y_i),(x_j,y_j)\big)\le \max\big\{\dist\big((x_0,y_0),(x_i,y_i)\big),\dist\big((x_0,x_j),(y_0,y_j)\big)\big\}.$$
Since the metric on $X$ has the property $\NP[1]$, there are two distinct numbers $i,j\in\{1,2,3\}$ such that $$d_X(x_i,x_j)\le\max\{d_X(x_0,x_i),d_X(x_0,x_j)\}.$$ The $\NP[0]$-property of the metric space $Y$ ensures that
$$d_Y(y_i,y_j)\le\max\{d_Y(y_0,y_i),d_Y(y_0,y_j)\}.$$Combining these two inequalities, we conclude that 
\begin{multline*}
\dist\big((x_i,y_i),(x_j,y_j)\big)=\max\{d_X(x_i,x_j),d_Y(y_i,y_j)\}\le\\
\le \max\{d_X(x_0,x_i),d_X(x_0,x_j),d_Y(y_0,y_i),d_Y(y_0,y_j)\}=\\
=\max\big\{d\big((x_0,y_0),(x_i,y_i)\big),d\big((x_0,x_j),(y_0,y_j)\big)\big\}.
\end{multline*}
\end{proof}

Lemma~\ref{P1+0} implies that for a positive real number $a$  the metric
$$\dist((x,y),(x',y'))=\max\{|x-x'|,|y-y'|\}$$on the product
$\II_a=[-1,1]\times\{0,a\}\subset\IR\times\IR$ has the Nagata property $\NP[1]$.

For a metric space $X$ we shall write $\II_a\hookrightarrow X$ if $X$ contains an isometric copy of the space $\II_a$.

\begin{lemma}\label{univ} For any $\GP[1]$-metric space $X$ the set $A=\{a\in(\frac1{16},\frac18):\II_a\hookrightarrow X\}$ has cardinality $|A|\le \dens(X)$.
\end{lemma}

\begin{proof}  For every $a\in A$ fix an isometric embedding $h_a:\II_a\to X$ and define a map $f_a:\II_1\to X$ by 
letting $f_a:(x,t)\mapsto h_a(x,at)$ for $(x,t)\in \II_1$. The map $f_a$ can be considered as an element of the function space $C(\II_1,X)$ endowed with the sup-metric $$\dist(f,g)=\sup_{t\in \II_1}\dist(f(t),g(t)).$$ By \cite[3.4.16]{En}, the density of  the function space $C(\II_1,X)$ is equal to the density of $X$. Now the assertion of the theorem will follow as soon as we check that the set $\F_A=\{f_a:a\in A\}$ is discrete in $C(\II_1,X)$. This will follow as soon as we show that $\dist(f_a,f_b)\ge \frac1{32}$ for any numbers $a\ne b$ in $A$.

To this end we first introduce some notation. For $a\in A$ and $i\in\{0,1\}$ let $$I_a^i=f_a([-1,1]\times\{i\}),\;\partial I_a^i=f_a(\{-1,1\}\times\{i\}),\;J_a^i=I_a^i\setminus\partial I_a^i,\;c_a^i=f_a(\textstyle{\frac12},i),$$ and $$V_a^i=\{x\in X:\dist(x,I_a^i)<\frac13\dist(x,\partial I_a^i)\}.$$ By Lemma~\ref{slice}, there is a non-expanding retraction $r^i_a:V^i_a\to J_a^i$ such that for every $x\in V_a^i$ and  $t\in J_a^i$ we get
\begin{equation}\label{formula}
\dist(x,t)=\max\{\dist(r^i_a(x),t),\dist(x,I_a^i)\}.\end{equation}
Moreover, for any points $x,y\in V_a^i$ with $\dist(x,I_a^i)\ne \dist(y,I_a^i)$ we get \begin{equation}\label{formula2}d(x,y)\ge \max\{\dist(x,I_a^i),\dist(y,I_a^i)\}.
\end{equation}

To derive a contradiction, assume that $\dist(f_a,f_b)<\e=\frac1{32}$ for some distinct numbers $a,b\in A$. Observe that
$$\dist(c_b^0,I^1_a)\le \dist(c_b^0,c_a^0)+\dist(c_a^0,I^1_a)<\e+a<\frac1{32}+\frac18=\frac{5}{32}$$while
$$\dist(c_b^0,\partial I^1_a)\ge \dist(c_a^0,\partial I^1_a)-\dist(c_a^0,c_b^0)=\frac12-\e=\frac12-\frac1{32}=\frac{15}{32}.$$
Consequently, $\dist(c_b^0,I^1_a)<\frac13\dist(c_b^0,\partial I^1_a)$ and hence $c_b^0\in V_a^1$. We claim that $\dist(c_b^0,I^1_a)=\dist(c_a^0,I^1_a)=a$. Otherwise, we may apply the formula~(\ref{formula2}) to derive a contradiction: $$d(c_b^0,c_a^0)\ge \max\{\dist(c_b^0,I_a^1),\dist(c_a^0,I_a^1)\}\ge\dist(c_a^0,I_a^1)=a>\e>\dist(f_a,f_b).$$  Since the retraction $r_a^1:V_a^1\to J_a^1$ is non-expanding, we get $$\dist(r_a^1(c_b^0),c^1_a)=\dist(r_a^1(c_b^0),r_a^1(c_a^0))\le\dist(c_b^0,c_a^0)<\e<a.$$ Now the formula~(\ref{formula}) yields
$$\dist(c_b^0,c_a^1)=\max\{\dist(r_a^1(c_b^0),c_a^1),\dist(c_b^0,I_a^1)\}=\dist(c_b^0,I_a^1)=a.$$

By analogy we can prove that $\dist(c_a^1,c_b^0)=b$, which contradicts $\dist(c_b^0,c_a^1)=a$.
\end{proof}

\section{Obtuse arcs and embeddings with small distortion}

In this section we shall introduce the notion of an obtuse arc and show that for each embedding $f:[0,1]\to X$ with $\Dist(f)<2$ the arc $f([0,1])$ is obtuse.
By a {\em metric arc} we understand a metric space that is homeomorphic to the unit interval $\II=[0,1]$.

A metric arc $(I,d)$ is called {\em obuse} if
\begin{itemize}
\item for any subarc $J\subset I$ with end-points $a,b$ and any point $z\in J\setminus\{a,b\}$ there are points $x,y\in J$  with $d(x,y)>\max\{d(z,x),d(z,y)\}$; and 
\item for any subarc $J\subset I$ with end-points $a,b$ there is a point $z\in J$ with $d(a,b)>\max\{d(z,a),d(z,b)\}$.
\end{itemize}
In this case the metric $d$ on $I$ is called {\em obtuse}.

It is easy to see that each subinterval $[a,b]\subset \IR$ endowed with the Euclidean metric is an obtuse arc. It can be shown that each continuously differentiable curve can be covered by finitely many obtuse subarcs.

\begin{proposition}\label{BM} If an embedding $f:\II\to X$ of the unit interval $\II=[0,1]$ into a metric space $(X,d_X)$ has distortion $\Dist(f)<2$, then the image $I=f(\II)$ is an obtuse arc in $X$.
\end{proposition}

\begin{proof}  We need to show that the metric $$\rho(t,t')=d_X(f(t),f(t'))$$ on $\II$, induced by the embedding $f$, is obtuse. It follows that $$(\|f^{-1}\|_\Lip)^{-1}\cdot|x-y|\le \rho(x,y)\le \|f\|_\Lip\cdot|x-y|.$$ Now we establish the two conditions of the definition of an obtuse arc.
\smallskip

1) Take any subinterval $[a,b]\subset\II$ and a point $z\in(a,b)$. Let $x,y\in(a,b)$ be any two points such that $z$ is the midpoint of the interval $(x,y)$. Then
\begin{multline*}
\max\{\rho(x,z),\rho(y,z)\}\le \|f\|_\Lip\cdot\max\{|x-z|,|y-z|\}=\|f\|_\Lip\cdot|x-y|/2\le\\
\le \frac12\|f\|_\Lip\cdot\|f^{-1}\|_\Lip\cdot\rho(x,y)<\frac12\cdot 2\cdot\rho(x,y)<\rho(x,y).
\end{multline*}
\smallskip

2) By analogy we can prove that for any subinterval $[a,b]\subset\II$ the midpoint $z$ of  $[a,b]$ satisfies the inequality $\max\{\rho(x,z),\rho(y,z)\}<\rho(x,y)$.
\end{proof}

\section{Obtuse arcs in $\NP[1]$-metric spaces}

In this section we study the structure of an $\NP[1]$-metric space $X$ in a neighborhood of an obtuse arc $I\subset X$.

\begin{proposition}\label{obtuse} Let $(X,d)$ be an $\NP[1]$-metric space, $I\subset X$ be an obtuse arc with endpoints $a,b$ in $X$ and let $V=\{x\in X:\dist (x,I)<\dist (x,\{a,b\})\}$.
\begin{enumerate}
\item For every point $x\in V\setminus I$  the set
$D(x)=\{t\in I:d(x,t)=d(x,I)\}$ is the finite union of closed subintervals of $I$ each of which has diameter  $>\dist(x,I)$.
\item For any points $x,y\in V$ with $\dist(x,I)\ne\dist(y,I)$ we get
$$d(x,y)\ge \max\{\dist(x,I),\dist(y,I)\}.$$
\end{enumerate}
\end{proposition}

\begin{proof} 1. Given a point $x\in V\setminus I$ put $D=\dist(x,I)$ and consider the family $\mathcal I$ of maximal non-generate subintervals in the closed subset $$D(x)=\{t\in I:d(t,x)=D\}\subset (a,b)=I\setminus\{a,b\}.$$

We claim that each maximal subinterval $[a_1,b_1]\in\mathcal I$ has diameter $\diam[a_1,b_1]> D$. Assuming conversely that $\diam([a_1,b_1])\le D$, and using the second condition of the definition of an obtuse metric, we can find a point $x_0\in (a_1,b_1)$ such that $D\ge d(a_1,b_1)>\max\{d(a_1,x_0),d(b_1,x_0)\}$.  The maximality of the subinterval $[a_1,b_1]\subset D(x)\subset (a,b)$ implies the existence of points $x_1\in (a,a_1)\setminus D(x)$ and $x_2\in (b_1,b)\setminus D(x)$ such that  $\max\{d(x_1,x_0),d(x_2,x_0)\}<\min\{D,d(x_1,x_2)\}$. Now we see that
the quadruple of points $x_0,x_1,x_2,x_3=x$ witnesses that the metric $d$ on $X$ fails to have the Nagata property $\NP[1]$ because
$$
\begin{aligned}
&d(x_1,x_2)>\max\{d(x_1,x_0),d(x_2,x_0)\},\\
&d(x_1,x_3)>D\ge \max\{d(x_0,x_1),d(x_0,x_3)\}\mbox{ and }\\
&d(x_2,x_3)>D\ge\max\{d(x_0,x_2),d(x_0,x_3)\}.
\end{aligned}$$

Taking into account that any two distinct maximal subintervals in the family $\mathcal I$ are disjoint and have diameter $>D$, we conclude that the family $\mathcal I$ is finite. It remains to show that $D(x)=\cup\mathcal I$. Assuming the converse, we could find a point $x_0\in D(x)\setminus\cup\mathcal I$ and  a neighborhood $(a_1,b_1)\subset I\setminus\cup\mathcal I$ of the point $x_0$ in $I\setminus\{a,b\}$ such that $\diam(a_1,b_1)<D$. The intersection $(a_1,b_1)\cap D(x)$ contains no non-degenerate subinterval and hence is nowhere dense in $(a_1,b_1)$. The obtuse property of the metric $d$ guarantees the existence of two points $x_1,x_2\in (a_1,b_1)$ such that $d(x_1,x_2)>\max\{d(x_1,x_0),d(x_2,x_0)\}$. Since $D(x)\cap (a_1,b_1)$ is nowhere dense we can additionally assume that $x_1,x_2\notin D(x)$. Then for the quadruple of the points $x_0,x_1,x_2,x_3=x$ we get
$$
\begin{aligned}
&d(x_1,x_2)>\max\{d(x_1,x_0),d(x_2,x_0)\},\\
&d(x_1,x_3)=d(x_1,x)>D=\max\{d(x_1,x_0),d(x_3,x_0)\}, \mbox{ and}\\
&d(x_2,x_3)=d(x_2,x)>D=\max\{d(x_3,x_0),d(x_2,x_0)\},
\end{aligned}
$$ witnessing the failure of the Nagata property $\NP[1]$ for the metric $d$.
\smallskip

2. Given two points $x,y\in V$ with $\dist(x,I)\ne\dist(y,I)$ we should prove that
$d(x,y)\ge\max\{\dist(x,I),\dist(y,I)\}$. Assume conversely, that
 $d(x,y)<\max\{d(x,I),d(y,I)\}$. Without loss of generality $d(x,I)<d(y,I)$. By the preceding item, the set $$D(x)=\{z\in I:d(x,z)=d(x,I)\}$$ contains two points $x_1,x_2$ with $d(x_1,x_2)>d(x,I)$.  Now we see that the quadruple of the points $x_0=x,x_1,x_2,x_3=y$ satisfies the inequalities
$$
\begin{aligned}
&d(x_1,x_2)>d(x,I)=\max\{d(x_0,x_1),d(x_0,x_2)\},\\
&d(x_1,x_3)\ge d(y,I)>\max\{d(x_0,x_1),d(x_0,x_3)\},\\
&d(x_2,x_3)\ge d(y,I)>\max\{d(x_0,x_2),d(x_0,x_3)\},
\end{aligned}
$$witnessing that the metric $d$ fails to have the Nagata property $\NP[1]$.
\end{proof}

By an argument similar to that from Theorem~\ref{c1}, we apply Proposition~\ref{obtuse} to prove the following

\begin{corollary}\label{c8} Let $X$ be a  locally connected $\NP[1]$-metric space $X$ and $I\subset X$ is an obtuse arc with endpoints $a,b$. Then the set $I\setminus\{a,b\}$ is open in $X$.
\end{corollary}

\section{Acknowledgment}

The authors express their sincere thanks to Michael Zarichnyi for fruitful discussions related to metric spaces with the Nagata property. We also acknowledge
the referee for comments and suggestions.


\begin{thebibliography}{99}

\bibitem{naor} Y.~Bartal, N.~Linial, M.~Mendel, A.~Naor, {\em On Some Low Distortion Metric Ramsey Problems}, Discrete Comput.~Geom.~{\bf 33}:1 (2005), 27--45.

\bibitem{Ellis} R.~Ellis, {\em Extending uniformly continuous pseudo-ultrametrics and uniform retracts}, Proc. Amer. Math. Soc. {\bf 30}:3 (1971), 599--602.

\bibitem{En} R.~Engelking, {\em General Topology}, Warsaw, PWN, 1977.

\bibitem{Eng} R.~Engelking, {\em Theory of Dimensions, Finite and Infinite}, Heldermann Verlag, 1995.

\bibitem{dG} J.~de Groot, {\em Non-archimedean metrics in topology}, Proc. Amer. Math. Soc. {\bf 7} (1956), 948--953.

\bibitem{LL} A.~Lemin, V.~Lemin, {\em On a universal ultrametric space}, Topology Appl. {\bf 103}:3 (2000), 339--345.

\bibitem{mat1} J.~Matou\v sek, {\em Ramsey-like properties for bi-Lipschitz mappings of finite metric spaces}, Comment. Math. Univ. Carolin. {\bf 33}:3 (1992), 451--463.

\bibitem{mat2} J.~Matou\v sek, {\em Note on bi-Lipschitz embeddings into normed spaces}, Comment. Math. Univ. Carolin. {\bf 33}:1 (1992), 51--55.

\bibitem{Men} M.~Mendel, {\em Metric Dichotomies}, in: Limits of graphs in group theory and computer science (G.~Arzhantseva, A.~Valette eds.), EPFL Press, Lausanne, 2009;   (arXiv:0710.1994).

\bibitem{Nag} J.~Nagata, {\em On a special metric and dimension}, Fund. Math. {\bf 55} (1964), 181--194.

\bibitem{NagB} J.~Nagata, {\em Modern Dimension Theory}, Sigma Series in Pure Math., 2. Heldermann Verlag, Berlin, 1983.

\bibitem{TZ} E.~D.~Tymchatyn, M.~Zarichnyi, {\em A note on operators extending partial ultrametrics}, Comment. Math. Univ. Carolin. {\bf 46}:3 (2005), 515--524.

\end{thebibliography}
\end{document}